\documentclass[a4paper,11pt]{amsart}
\usepackage{amssymb}
\usepackage{amsmath}
\usepackage{mathrsfs}
\usepackage{amsmath,amssymb,amsthm,latexsym,amscd,mathrsfs}
\usepackage{indentfirst}
\usepackage{graphicx}
\usepackage{booktabs}
\usepackage{array}
\usepackage{chemarrow}
\newcommand{\ROM}[1]{\mathrm{\uppercase\expandafter{\romannumeral#1}}}
\theoremstyle{definition}

\newtheorem{thm}{Theorem}[section]
\newtheorem{lem}{Lemma}[section]

\newtheorem{rem}{Remark}[section]

\makeatletter


\title[Isoparametric foliations and critical sets of eigenfunctions]{\textbf{Isoparametric foliations and critical sets of eigenfunctions}}
\author[Z.Z.Tang]{Zizhou Tang}\address{School of Mathematical Sciences, Laboratory of Mathematics and Complex Systems, Beijing Normal
University, Beijing 100875, China}\email{zztang@bnu.edu.cn}
\thanks {The project is partially supported by the NSFC ( No.11331002, No.11301027) and SRFDP (No. 20130003120008).}
\author[W. J. Yan]{Wenjiao Yan$^{\dag}$}
\address{School of Mathematical Sciences, Laboratory of Mathematics and Complex Systems, Beijing Normal
University, Beijing 100875, China} \email{wjyan@bnu.edu.cn}
\thanks {$^{\dag}$ the corresponding author}

\subjclass[2010]{ 53C40, 57R70, 58J50.}
\date{}
\keywords{eigenfunction of Laplacian, isoparametric hypersurface,
critical set.}
\begin{document}

\maketitle
\begin{abstract}
Jakobson and Nadirashvili \cite{JN} constructed a sequence of eigenfunctions
on $T^2$ with a bounded number of critical points, answering in the negative
the question raised by Yau \cite{Yau1} which asks that whether the number of the critical points
of eigenfunctions for the Laplacian increases with the corresponding
eigenvalues.

The present paper finds three interesting eigenfunctions
on the minimal isoparametric hypersurface $M^n$ in
$S^{n+1}(1)$. The corresponding
eigenvalues are $n$, $2n$ and $3n$, while their critical sets consist of
$8$ points, a submanifold(infinite many points) and $8$ points, respectively. On one of
its focal submanifolds, a similar phenomenon occurs.

\end{abstract}

\section{Introduction}

Eigenvalues of Laplacian are very important intrinsic invariants,
which reflect the geometry of manifolds very precisely.
Unfortunately, there are few manifolds whose eigenvalues are clearly
known, not to mention the eigenfunctions. The numbers of critical
points of eigenfunctions are even more difficult to determine.
However, as S.T.Yau pointed out, this number is closely related to
many important questions, which makes it worthy of being studied
extensively. In this regard, S.T.Yau \cite{Yau1} raised a question:
is it true that the number of critical points of the $k$-th eigenfunction on a
compact Riemannian manifold increases with $k$. He also investigated this problem in
the surface case (cf. \cite{Yau2}).

In 1999, Jakobson and Nadirashvili \cite{JN} constructed
a metric on a $2$-dimensional torus and a sequence of eigenfunctions
such that the corresponding eigenvalues go to infinity while the
number of critical points remains bounded, a constant in fact. But
in a fastidious manner, this remarkable example does not deny Yau's conjecture
in the sense of ``non-decreasing".

In the present paper, by taking advantage of a natural
concept--isoparametric hypersurface, we find (based on \cite {Sol})
an isoparametric function, which is an eigenfunction on the minimal
isoparametric hypersurface $M^n$ of OT-FKM type in $S^{n+1}(1)$.
Combining with the other two well-known eigenfunctions, it
constitutes a sequence of eigenfunctions with increasing eigenvalues,
but the numbers of their critical points are not monotonic at all.

Similarly, another isoparametric function (indeed an eigenfunction)
expressed in the same form arises in one of the focal submanifolds
of $M^n$ mentioned before. Together with the other eigenfunction,
it constitutes a sequence of eigenfunctions with similar property
as that on $M^n$.

One of the main results of the present paper is the following:
\vspace{1mm}

\begin{thm}\label{thm1}
\emph{Let $M^n$ be the minimal isoparametric hypersurface of
OT-FKM type in the unit sphere $S^{n+1}(1)$. Then there exist three
eigenfunctions $\varphi_1$, $\varphi_2$ and $\varphi_3$ defined on
$M^n$, corresponding to eigenvalues $n$, $2n$ and $3n$, whose
critical sets consist of $8$ points, a submanifold and $8$ points,
respectively. For specific, $\varphi_1$ and $\varphi_3$ are both
Morse functions; $\varphi_2$ is an isoparametric function on $M^n$,
whose critical set $C(\varphi_2)$ is:
\begin{equation}\label{varphi2}
C(\varphi_2)=N_+\cup N_-,\quad \dim N_+=\dim N_-= n-m ~(1\leq m <
n),
\end{equation}
where the number $m$ will be introduced in the definition of
OT-FKM type.
}
\end{thm}

\begin{rem}
The Morse number (the minimal number of critical points of all Morse
functions) of a compact isoparametric hypersurface with $g=4$
distinct principal curvatures in the unit sphere is equal to $2g=8$
(\emph{cf.} \cite{CR}).
\end{rem}
\vspace{1mm}

Firstly, to clarify notations, we denote the Laplacian on an
$n$-dimensional compact manifold $M^n$ by $\Delta f= \text{div}
\nabla f$, and say $\lambda_k$ its $k$-th eigenvalue with
multiplicity ($\lambda_0=0<\lambda_1<\lambda_2<...$) if $\Delta f_k
+ \lambda_k f_k =0$ for some $f_k: M^n\rightarrow \mathbb{R}$.
Correspondingly, $f_k$ is called the $k$-th eigenfunction. The
present paper is mainly concerned with the number of critical points
of the eigenfunction $f_k$.

Recall that a hypersurface $M^n$ in a Riemannian manifold
$\widetilde{M}^{n+1}$ is \emph{isoparametric} if it is a level
hypersurface of an isoparametric function $f$ on
$\widetilde{M}^{n+1}$, that is, a non-constant smooth function $f:
\widetilde{M}^{n+1}\rightarrow \mathbb{R}$ satisfying (cf. \cite{Wan, {GT2}}):
\begin{equation}\label{ab}
\left\{ \begin{array}{ll}
~~|\widetilde{\nabla} f|^2= b(f)\\
\quad\widetilde{\triangle} f~~=a(f)
\end{array}\right.
\end{equation}
where $b$ and $a$ are smooth and continuous functions on $\mathbb{R}$, respectively.

In this meaning, 
the \emph{focal varieties} are the preimages of the global maximum
and minimum values (if exist) of $f$, which we denote by $M_+$ and $M_-$,
respectively. They are in fact both minimal submanifolds of
$\widetilde{M}^{n+1}$ with codimensions $m_++1$ and $m_-+1$ in
$\widetilde{M}^{n+1}$, respectively(cf. \cite{Wan},\cite{Th},
\cite{GT1}).

As asserted by \'{E}lie Cartan, an isoparametric hypersurface in the
unit sphere is indeed a hypersurface with constant principal
curvatures. Let $g$ be the number of distinct principal curvatures with multiplicity $m_i$ ($i=1, \cdots, g$).
An elegant result of M\"{u}nzner states that $g$ can be only $1, 2, 3, 4$ or $6$, and $m_i=m_{i+2}$ (subscripts mod $g$).
To clarify the notations, we denote $m_+=:m_1$ and $m_-=:m_2$.
Up to now, the isoparametric hypersurfaces with $g = 1, 2, 3, 6$ are
completely classified (cf. \cite{DN} and \cite{Miy}). For
isoparametric hypersurfaces with $g=4$, Cecil-Chi-Jensen
(\cite{CCJ}), Immervoll (\cite{Imm}) and Chi (\cite{Chi, Chi2}) proved a
far reaching result that they are all of OT-FKM type except for the homogeneous case
with $(m_+, m_-)=(2, 2), (4, 5)$.

From now on, we are specifically concerned with the isoparametric hypersurfaces of OT-FKM type in $S^{n+1}(1)$
with four distinct principal curvatures. For a symmetric Clifford system $\{P_0,...,P_m\}$
on $\mathbb{R}^{2l}$, \emph{i.e.}, $P_i$'s are symmetric matrices
satisfying $P_iP_j+P_jP_i=2\delta_{ij}I_{2l}$, the \emph{OT-FKM type
isoparametric hypersurfaces} are level hypersurfaces of $f:=F|_{S^{2l-1}}$ with $F$ defined
by  Ferus, Karcher and M\"{u}nzner (cf. \cite{FKM}):
\begin{eqnarray}\label{OT-FKM isop. poly.}
&& \quad F: \mathbb{R}^{2l} \rightarrow \mathbb{R}\nonumber\\
&& F(x) = |x|^4 - 2\displaystyle\sum_{\alpha = 0}^{m}{\langle
P_{\alpha}x,x\rangle^2}
\end{eqnarray}
The pairs $(m_+, m_-)$ of the OT-FKM type are $(m, l-m-1)$, provided
$m>0$ and $l-m-1>0$, where $l = k\delta(m)$ $(k=1,2,3,...)$, $\delta(m)$ is the dimension
of an irreducible module of the Clifford algebra $C_{m-1}$, which we list below:
\vspace{1mm}

\begin{center}
\begin{tabular}{|c|c|c|c|c|c|c|c|c|c|}
\hline
$m$ & 1 & 2 & 3 & 4 & 5 & 6 & 7 & 8 & $\cdots$ $m$+8 \\
\hline
$\delta(m)$ & 1 & 2 & 4 & 4 & 8 & 8 & 8 & 8 & ~16$\delta(m)$\\
\hline
\end{tabular}
\end{center}
\vspace{1mm}

We now fix $M^n$ to be the minimal isoparametric hypersurface of
OT-FKM type in $S^{n+1}(1)$, and $f$ to be $f:=F|_{S^{2l-1}}$ with $F$
defined in (\ref{OT-FKM isop. poly.}). Choosing a point $q_1\in
S^{n+1}(1)\backslash \{M_+, M_-, M^n\}$, we define three
eigenfunctions $\varphi_1$, $\varphi_2$ (following \cite{Sol}) and
$\varphi_3$ as follows:
\begin{eqnarray}
&& \varphi_1: M^n\rightarrow \mathbb{R},\qquad\qquad \varphi_2: M^n\rightarrow \mathbb{R}\qquad\qquad \varphi_3: M^n\rightarrow \mathbb{R}\nonumber\\
&& \qquad x \mapsto \langle x, q_1\rangle, \qquad\qquad\quad x \mapsto \langle Px, x\rangle \quad\qquad\qquad x \mapsto \langle \xi(x), q_1\rangle
\end{eqnarray}
\vspace{1mm} where $\xi$ is a unit normal vector field on $M^n$;
$P\in \Sigma:=\Sigma (P_0,...,P_m)$, the unit sphere in
$\mathrm{\mathrm{Span}}\{P_0,...,P_m\}$, which is called \emph{the Clifford sphere}
(see Definition 3.6 of \cite{FKM}).

\noindent
\begin{rem}
It was proved by the authors that the first eigenvalue of the closed
minimal isoparametric hypersurface $M^n$ in $S^{n+1}(1)$ is just $n$
(cf. \cite{TY}, \cite{TXY}). As a corollary, the coordinate function
restricted on $M^n$, $\varphi_1$, is the first eigenfunction.
\end{rem}

With all the preconditions, a direct verification reveals that the
eigenvalues corresponding to $\varphi_1$, $\varphi_2$ and
$\varphi_3$ are $n$, $2n$ and $3n$, respectively.
Moreover, with our choice of $q_1\in S^{n+1}(1)\backslash \{M_+,
M_-, M^n\}$, a simple application of isoparametric geometry shows
that $\varphi_1$ and $\varphi_3$ are both Morse functions with
$2g=8$ critical points. The more fascinating result is that
$\varphi_2$ is indeed an isoparametric function on $M^n$, thus by
virtue of \cite{Wan}, the critical set of $\varphi_2$ are just the
union of its focal submanifolds $N_+$ and $N_-$. For the proof of
Theorem \ref{thm1}, we need the following lemma

\begin{lem}\label{N+N-}
\emph{For focal submanifolds $N_+$ and $N_-$ of $\varphi_2$ on
$M^n$, we have diffeomorphisms:
\begin{equation*}
N_+\underset{diff.}{\cong}N_-\underset{diff.}{\cong}M_+ =\{x\in
S^{n+1}(1)~|~\langle P_0x, x\rangle=\langle P_1x,
x\rangle=\cdots=\langle P_mx, x\rangle=0\}.
\end{equation*}
Particularly, in the case of $m=1$, each level
(isoparametric) hypersurface of $\varphi_2$ is minimal in $M^n$.}
\end{lem}

\begin{rem}
When $m=1$, the codimensions of $N_+$ and $N_-$ in $M^n$ are $1$, this is what called improper isoparametric (cf. pp.165 of \cite{GT2}).
\end{rem}

As we stated before, another counterexample of Yau's
conjecture appear on the focal submanifold $M_-:=f^{-1}(-1)$ with
dimension $l+m-1$. In a similar way, we define two eigenfunctions
$\omega_1$, $\omega_2$ on $M_-$:
\begin{eqnarray}\label{omega12}
&& \omega_1: M_- \rightarrow \mathbb{R}\qquad\qquad\qquad\qquad \omega_2: M_-\rightarrow \mathbb{R}\nonumber\\
&& \qquad x \mapsto \langle Px, x\rangle\qquad\qquad\quad\qquad\qquad x \mapsto \langle x, q_2\rangle,
\end{eqnarray}
where $P\in \Sigma=\Sigma(P_0, P_1,...,P_m)$, $q_2\in
S^{n+1}(1)\backslash \{M_+, M_-\}$. Correspondingly, we have the
following theorem:

\begin{thm}\label{thm2}
\emph{Let $M_-:=f^{-1}(-1)$, a focal submanifold of OT-FKM type in
the unit sphere $S^{n+1}(1)$. Then there exist two eigenfunctions
$\omega_1$ and $\omega_2$ defined on $M_-$, corresponding to
eigenvalues $4m$ and $l+m-1$, whose critical sets consist of a
submanifold and $4$ points, respectively.  For specific, $\omega_2$
is a Morse function; $\omega_1$ is an isoparametric function on
$M_-$, whose critical set $C(\omega_1)$ is:
\begin{equation}\label{omega1}
C(\omega_1)=V_+\cup V_-, \quad \dim V_+=\dim V_-= l-1.
\end{equation}}
\end{thm}

\begin{rem}
The Morse number of each focal submanifold of a compact
isoparametric hypersurface with $g=4$ distinct principal curvatures
in the unit sphere is equal to $g=4$ (\emph{cf.} \cite{CR}).
\end{rem}
\vspace{1mm}

For the proof of Theorem \ref{thm2}, we need the following:

\begin{lem}\label{V+V-}
\emph{For focal submanifolds $V_+$ and $V_-$ of $\omega_1$ on $M_-$,
we have isometries:
\begin{equation*}
V_+\underset{isom.}{\cong}V_-\underset{isom.}{\cong}S^{l-1}(1).
\end{equation*}
Particularly, in the improper case, i.e. $m=1$, each level (
isoparametric ) hypersurface of $\omega_1$ is minimal in $M_-$.}
\end{lem}

Comparing with the values of $\delta(m)$ in the previous table, we
observe that $4m<l+m-1$ at most cases. More precisely, $4m<l+m-1$ as
long as $k\geq 5$ and $m\leq 9$; $4m<l+m-1$ holds true for any $k$
when $m\geq 10$. Therefore, with an appropriate choice of $k$, we
can always make eigenfunctions $\omega_1$ and $\omega_2$ another
counterexample of Yau's conjecture.

Bearing these examples in mind, we would like to raise the following
question:
\vspace{1mm}

\noindent \textbf{Question:} For a generic metric on a compact
manifold $M$, is the number of critical points of the first
eigenfunction (must be a Morse function, according to Uhlenbeck \cite{Uh})
equal to the Morse number of $M$?\footnote{\emph{Added in proof}. It was recently proved by A. Enciso and D. Peralta-Salas that on a compact manifold, there is a Riemannian metric such that the first nontrivial eigenfunction can have as many non-degenerate critical points as one wishes (bigger in particular than the Morse number of the manifold). Moreover, any other metric $C^{\infty}$ close to it carries the same property(cf. \cite{EP-S}).}

\section{counterexamples on $M^n$}

This section will be committed to proving Theorem \ref{thm1} on the
minimal isoparametric hypersurface $M^n$ of OT-FKM type in
$S^{n+1}(1)$. At first, we denote the connections and Laplacians on
$M^n$, $S^{n+1}(1)$ and $\mathbb{R}^{n+2}$ respectively by:
\begin{eqnarray*}\label{connection}
&& M^n \subset S^{n+1}(1) \subset \mathbb{R}^{n+2}\nonumber\\
&& \nabla~ \triangle, \quad\overline{\nabla}~ \overline{\triangle}, \qquad\widetilde{\nabla} ~\widetilde{\triangle}.
\end{eqnarray*}

In order to facilitate the description, we state the following lemma in front of the proof of Theorem \ref{thm1}. The proof
is direct and will be omitted here.
\begin{lem}\label{laplacian relation}
\emph{Let $\xi$ be a ( local )unit vector field on $S^{n+1}(1)$
extended from a unit normal vector field of $M^n$, $H$ be the mean
curvature vector field of $M^n$ in $S^{n+1}(1)$. For functions $\mathcal{G}$ on
$\mathbb{R}^{n+2}$, $G=\mathcal{G}|_{S^{n+1}}$ and $g=G|_{M^n}$, at
any $x\in M^{n}$ ( as a position vector field ) we have:
\begin{equation}\label{two laplacian}
\left\{ \begin{array}{ll}
\widetilde{\triangle}\mathcal{G}|_{S^{n+1}} = \overline{\triangle} G+nx(\mathcal{G})+xx(\mathcal{G})\\
\overline{\triangle}G|_{M^n} ~~= \triangle g -\xi(G)\langle H,
\xi\rangle + \xi\xi(G) -\overline{\nabla}_{\xi}\xi(G)
\end{array}\right.
\end{equation}}
\end{lem}

\vspace{4mm}

\noindent \textbf{\emph{Proof of Theorem \ref{thm1}.}}\quad We take
the first step by determining the eigenvalues corresponding to
$\varphi_i$ ($i=1, 2, 3$).
Clearly, based on Lemma \ref{laplacian relation}, a direct calculation depending on the
minimality of $M^n$ in $S^{n+1}(1)$ leads to
\begin{equation}\label{eigenvalue of varphi13}
\triangle \varphi_1= -n \varphi_1.
\end{equation}
Besides, in conjunction with Codazzi equation, we get another straightforward result:
\begin{equation}\label{eigenvalue of varphi31}
\triangle \varphi_3= -|B|^2 \varphi_3 = -(g-1)n \varphi_3 = -3n \varphi_3,
\end{equation}
where $B$ is the second fundamental form of $M^n$, and the second equality in (\ref{eigenvalue of varphi31}) is an assertion
of \cite{PT}.
According to Solomon \cite{Sol}, the eigenvalue corresponding to
$\varphi_2$ is equal to $2n$. As a matter of fact, this conclusion
can also be derived from a few basic facts and Lemma \ref{laplacian
relation}---some formulas in this process will be useful later:

It is well known that there exists a unique $c_0$ with $-1<c_0<1$
such that the minimal isoparametric hypersurface $M^n$ (of OT-FKM type)
is given by $M^n=f^{-1}(c_0)$ (the value of $c_0$ will be given in the proof of Lemma \ref{N+N-}). We can choose the unit normal vector
field to be
$$\xi=\frac{\overline{\nabla}f}{|\overline{\nabla}f|}\Big|_{M^n}=\frac{\widetilde{\nabla}F-4Fx}{4\sqrt{1-F^2}}\Big|_{M^n}.$$
Extending $\xi$ along the normal geodesics such that
$\overline{\nabla}_{\xi}\xi=0$, it follows that
\begin{equation}\label{xi varphi2} \xi(\varphi_2) = \langle \xi,
\overline{\nabla}\varphi_2 \rangle =
\langle\frac{\widetilde{\nabla}F-4Fx}{4\sqrt{1-F^2}},
2Px-2\varphi_2x \rangle = -2\sqrt{\frac{1+f}{1-f}}\varphi_2,
\end{equation}
and thus
\begin{equation*}
\xi\xi(\varphi_2) = \langle \xi,
\overline{\nabla}\xi(\varphi_2)\rangle = -4\varphi_2.
\end{equation*}
Here, we extended $\varphi_2$ to $S^{n+1}(1)$ and $\mathbb{R}^{n+2}$
in a natural way. Then combining with (\ref{two laplacian}) and
$H=0$, we arrive at
\begin{equation}\label{eigenvalue of varphi2}
\triangle \varphi_2= -2n \varphi_2.
\end{equation}

\vspace{4mm}

Next, we aim to investigate the critical sets of $\varphi_i$ ($i=1,
2, 3$). Let $e_1, e_2,...,e_n$ be an orthonormal tangent frame field
on $M^n$ with $A_{\xi}e_i=\mu_ie_i$ ($i=1,2,...,n$), where $A_{\xi}$
is the shape operator. According to M\"{u}nzner, the principal
curvature $\mu_i\in\{\cot \theta_j =
\cot(\theta_1+\frac{j-1}{4}\pi)~|~ 0<\theta_1<\frac{\pi}{4},~
j=1,2,3,4\}.$
\vspace{1mm}

\emph{\emph{(i)}} 
For each $e_i\in T_xM^n$, we have
\begin{equation}\label{gradient varphi1}
\langle \nabla \varphi_1, e_i\rangle = e_i\langle x, q_1\rangle =
\langle e_i, q_1\rangle.
\end{equation}
It follows that $x$~is~a~critical~point~of~$\varphi_1$ if and only
if $q_1\in \mathrm{Span}\{ x, \xi(x) \}$.
In other words, $q_1$ lies on some normal geodesic $v(t)$ ($-\pi\leq
t\leq \pi$) with $v(0)=x,~v^{\prime}(0)=\xi(x)$. Therefore the
number of critical points of $\varphi_1$ is
\begin{equation*}
\sharp C(\varphi_1)=\frac{2\pi}{\pi/g}=2g=8.
\end{equation*}
Here, we used the known fact that the distance between two focal
submanifolds is equal to $\pi/g$ (cf. \cite{CR}).
Furthermore, recall the formula of Hessian:
$$\mathrm{Hess}(\varphi_1)_{ij}=\langle e_i, \nabla_{e_j}\nabla\varphi_1\rangle.$$
Restricted to a critical point $x$, using (\ref{gradient varphi1})
we express it as
\begin{equation}\label{Hessian 1}
\mathrm{Hess} (\varphi_1)|_x=-\mathrm{diag}\{~\langle\mu_1\xi-x, q_1\rangle, ~\langle\mu_2\xi-x, q_1\rangle,...,~\langle\mu_n\xi-x, q_1\rangle~\}.
\end{equation}
Writing $q_1=\cos t~x+\sin t~\xi$ $(-\pi<t<\pi)$ for a fixed $x$, a
direct calculation leads to
\begin{eqnarray*}
\langle\mu_i\xi-x, q_1\rangle=0 &\Leftrightarrow& \sin t(\cot\theta_i-\cot t)=0\\
  &\Leftrightarrow& q_1\in M_+\cup M_-\cup M^n.
\end{eqnarray*}
From the assumption $q_1\in S^{n+1}(1)\backslash\{M_+, M_-, M^n\}$,
we derive that $\varphi_1$ is a Morse function, as desired.

\vspace{1mm}

\emph{\emph{(ii)}} Similarly, for each $e_i\in T_xM^n$, we have
\begin{equation*}
\langle \nabla \varphi_3, e_i\rangle = e_i\langle \xi, q_1\rangle = -\langle A_{\xi}e_i, q_1\rangle = -\langle \mu_ie_i, q_1\rangle.
\end{equation*}
Since $\mu_i\in\{\cot \theta_j = \cot(\theta_1+\frac{j-1}{4}\pi)~|~
0<\theta_1<\frac{\pi}{4},~ j=1,2,3,4\},$ it is easy to see that
$\mu_i\neq 0~ ~\forall i$. Thus
$x$~is~a~critical~point~of~$\varphi_3$ if and only if 
$q_1\in \mathrm{Span}\{x, \xi(x) \}$.
Analogously, \begin{equation*}
\sharp C(\varphi_3)=\frac{2\pi}{\pi/g}=2g=8.
\end{equation*}
Furthermore, $\mathrm{Hess}(\varphi_3)$ at a critical point $x$ can be
expressed as
\begin{equation}\label{Hessian 3}
\mathrm{Hess}(\varphi_3)|_x=-\mathrm{diag}\{~\mu_1\langle\mu_1\xi-x, q_1\rangle, ~\mu_2\langle\mu_2\xi-x, q_1\rangle,...,~\mu_n\langle\mu_n\xi-x, q_1\rangle~\}.
\end{equation}
Again, our choice of $q_1$ guarantees that $\varphi_3$ is a Morse
function.

\vspace{1mm}

\emph{\emph{(iii)}} From the formula (\ref{xi varphi2}), we derive that
\begin{eqnarray}\label{nabla varphi2}
\nabla\varphi_2 &=& \widetilde{\nabla}\varphi_2-x(\varphi_2)x -\xi(\varphi_2)\xi \\
&=& 2(Px-\varphi_2x+\varphi_2\sqrt{\frac{1+c_0}{1-c_0}}\xi).\nonumber
\end{eqnarray}
Immediately, a simple calculation shows that $\varphi_2$ satisfies
\begin{equation}\label{isoparametric varphi2}
\left\{ \begin{array}{ll}
|\nabla \varphi_2|^2=4(1-\frac{2}{1-c_0}\varphi_2^2)\\
\quad\triangle \varphi_2 = -2n \varphi_2.
\end{array}\right.
\end{equation}
By definition, $\varphi_2$ is an isoparametric function on $M^n$.
Define the focal submanifolds by $N_{\pm}:=\{x\in
M^n~|~\varphi_2=\pm \sqrt{\frac{1-c_0}{2}}\}$. Therefore the
critical set of $\varphi_2$ is the union of its focal submanifolds:
$$C(\varphi_2)=N_+\cup N_-.$$

We are now in a position to complete the proof of Theorem \ref{thm1} by verifying Lemma \ref{N+N-}.

\vspace{2mm}

\noindent
\textbf{\emph{Proof of Lemma \ref{N+N-}.}}~
As indicated before, the focal submanifold $M_+$ of OT-FKM type is
$$M_+:=f^{-1}(+1)=\{x\in S^{n+1}(1)~|~\langle P_0x, x\rangle=\langle P_1x, x\rangle=\cdots=\langle P_mx, x\rangle=0\}.$$
Define a map:
\begin{eqnarray*}
&&h_+: M_+\rightarrow S^{n+1}(1) \\
&&\qquad\quad x\mapsto \cos t~ x+\sin t~ Px
\end{eqnarray*}
where $\cos t=\sqrt{\frac{1}{2}(1+\sqrt{\frac{1+c_0}{2}})}$, $\sin
t=\sqrt{\frac{1}{2}(1-\sqrt{\frac{1+c_0}{2}})}$. It is easy to show
that $$\langle Ph_+(x), h_+(x)\rangle=\sqrt{\frac{1-c_0}{2}},~
\emph{i.e.} ~~~~h_+(x)\in N_+.$$ Thus the image of $h_+$ is
contained in $N_+$. On the other hand, define another map:
\begin{eqnarray*}
&&j_+: N_+\rightarrow M_+ \\
&&\qquad\quad x\mapsto \cos t~ x+\sin t~ \xi(x)
\end{eqnarray*}
with the same values of $\cos t$ and $\sin t$, and
$\xi=\frac{\overline{\nabla}f}{|\overline{\nabla}f|}.$ Evidently,
$j_+$ is well defined and is just the inverse function of $h_+$.
This means that the focal submanifold $N_+$ of $\varphi_2$ on $M^n$
is diffeomorphic to the focal submanifold $M_+$ of $f$ on
$S^{n+1}(1)$.

We conclude the proof by investigating the mean curvatures of the
level hypersurfaces $N_t:=\varphi_2^{-1}(t)$, $t\in
(-\sqrt{\frac{1-c_0}{2}}, \sqrt{\frac{1-c_0}{2}})$. Following the
formula of the mean curvature $h(t)$ (cf. \cite{GT2}), we
have:
\begin{equation}\label{mean curvature}
h(t) = \frac{b^{\prime}(t)-2a(t)}{2\sqrt{b(t)}}
     = \frac{n-\frac{4}{1-c_0}}{\sqrt{1-\frac{2t^2}{1-c_0}}}~t
\end{equation}
Obviously, the isoparametric hypersurface $N_0=\varphi_2^{-1}(0)$ is
minimal in $M^n$. In addition, the minimality of $M^n$ implies:
$$c_0=\frac{m_--m_+}{m_-+m_+}=\frac{l-2m-1}{l-1},\quad n=2l-2,$$
then we obtain that
$$n-\frac{4}{1-c_0}=0 \Leftrightarrow m=1 ~(~the~ improper~ case~ (\mathrm{cf}.~ [GT2])).$$
In conclusion, in the improper case, all the level hypersurfaces of $\varphi_2$ are minimal.

The same argument applies to $N_-$ with a little change of the
values:
$$\cos t=\sqrt{\frac{1}{2}(1+\sqrt{\frac{1+c_0}{2}})}, \quad \sin t=-\sqrt{\frac{1}{2}(1-\sqrt{\frac{1+c_0}{2}})}.$$

 \hfill $\Box$

The proof of Theorem \ref{thm1} is now complete.

\section{counterexamples on $M_-$}
\noindent \textbf{\emph{Proof of Theorem \ref{thm2}.}} Implementing
the previous arguments in Section $2$, it is not difficult to find
that $\omega_2$ on $M_-$ is an eigenfunction corresponding to the
eigenvalue $\dim M_-=l+m-1$, and the number of its critical points
is $\frac{2\pi}{2\pi/g}=g=4$ (cf. \cite{CR}). Therefore, in
order to complete the proof of Theorem \ref{thm2}, we need only to
confirm that $\omega_1$ is an isoparametric function on $M_-$ and
prove Lemma \ref{V+V-}.

Firstly, noticing the Euclidean gradient $\widetilde{\nabla}\omega_1$ can be expressed by
$$\widetilde{\nabla}\omega_1=2Px=2\langle Px, x\rangle x+2\Big(Px-\langle Px, x\rangle x\Big),$$
we claim that
\vspace{1mm}

\noindent
\textbf{Claim}: \emph{$y:=Px-\langle Px, x\rangle x\in T_xM_-.$}
\vspace{1mm}

Holding this claim, it follows that $\nabla \omega_1=2y=2(Px-\langle
Px, x\rangle x)$. Then a simple calculation leads to
\begin{equation}\label{isoparametric omega1}
\left\{\begin{array}{ll}
|\nabla \omega_1|^2=4(1-\omega_1^2) \\
\quad\triangle \omega_1=-4m\omega_1,
\end{array}\right.
\end{equation}
where the second equality is due to Solomon \cite{Sol}. Namely,
$\omega_1$ is an isoparametric function on $M_-$. Define the focal
submanifolds of $\omega_1$ by $V_{\pm}:=\{x\in M_-~|~\omega_1=\pm
1\}$. Then the critical set of $\omega_1$ is
$$C(\omega_1)=V_+\cup V_-.$$

\begin{rem}
The proof of $|\nabla \omega_1|^2=4(1-\omega_1^2)$ is recently used by \cite{QT} to 
obtain a sequence of isoparametric functions (hypersurfaces).
\end{rem}
Now we are left to prove the previous Claim and Lemma \ref{V+V-}.

\vspace{1mm}

\noindent \textbf{\emph{Proof of \textbf{Claim}}}. Firstly, we
rewrite the focal submanifold
$$M_- := \{x\in S^{n+1}(1)~|~\displaystyle\sum_{\alpha=0}^{m}\langle P_{\alpha}x, x\rangle^2=1\}$$
as
\begin{eqnarray*}
M_- &=& \{x\in S^{n+1}(1)~|~x=\sum_{\alpha=0}^{m}\langle P_{\alpha}x, x\rangle P_{\alpha}x\}
\end{eqnarray*}
Define $\mathcal{P}:=\displaystyle\sum_{\alpha=0}^{m}\langle
P_{\alpha}x, x\rangle P_{\alpha}$, then for each $x\in M_-$ we have
\begin{equation}\label{Px=x}
\mathcal{P}\in \Sigma\quad and \quad\mathcal{P}x=x.
\end{equation}
Since $\mathcal{P}$ is an orthogonal symmetric matrix with vanishing trace,
we can decompose $\mathbb{R}^{2l}$ as
$$\mathbb{R}^{2l}=E_+(\mathcal{P})\oplus E_-(\mathcal{P}).$$
With respect to this decomposition, $2y\in\mathbb{R}^{2l}$ can be
written as
$$2y=(y+\mathcal{P} y)+(y-\mathcal{P}y).$$
Denoting $P=\displaystyle\sum_{\beta=0}^{m}a_{\beta}P_{\beta}$ with $\displaystyle\sum_{\beta=0}^{m}a_{\beta}^2=1$, we have
\begin{eqnarray*}
y+\mathcal{P} y &=& Px-\langle Px, x\rangle x + \mathcal{P}Px-\langle Px, x\rangle \mathcal{P}x\\
&=& P\mathcal{P}x  + \mathcal{P}Px - 2\langle Px, x\rangle x\\
&=& \sum_{\beta=0}^{m}a_{\beta}P_{\beta}\Big(\sum_{\alpha=0}^{m}\langle P_{\alpha}x, x\rangle P_{\alpha}x\Big)+ \sum_{\alpha=0}^{m}\langle P_{\alpha}x, x\rangle P_{\alpha}\Big(\sum_{\beta=0}^{m}a_{\beta}P_{\beta}x\Big)- 2\langle Px, x\rangle x\\
&=& 2\sum_{\alpha=0}^{m}a_{\alpha}\langle P_{\alpha}x, x\rangle x-2 \sum_{\beta=0}^{m}a_{\beta}\langle P_{\beta}x, x \rangle x\\
&=& 0,
\end{eqnarray*}
which leaves $2y=y-\mathcal{P}y$, \emph{i.e.} $y\in E_-(\mathcal{P})$.

On the other hand, setting $y=Px-\langle Px, x\rangle x=Qx$, where
$$Q:=P-\langle Px, x\rangle\mathcal{P}\in \mathrm{Span}\{P_0, P_1,...,P_m\},$$
it is easy to find that
$$\langle Q, \mathcal{P}\rangle=0.$$
Comparing with (cf. Section 4.5(iii) of \cite{FKM})
\begin{equation*}
T_x^{\perp}M_-=\{\nu\in E_-(\mathcal{P})~|~\langle \nu, Qx\rangle=0,
~\forall~ \langle Q, \mathcal{P}\rangle=0\},
\end{equation*}
we get immediately the Claim.\hfill$\Box$

Now we are in a position to prove Lemma \ref{V+V-}.

\noindent
\textbf{\emph{Proof of Lemma \ref{V+V-}.}}
Under an orthogonal transformation, we can express $P$ as
\begin{equation*}
    P=T^t\left(
           \begin{array}{cc}
             I_l & 0 \\
             0 & -I_l \\
           \end{array}
         \right)T,\quad with ~T^tT=I_{2l}.
\end{equation*}
Write $Tx=(z, w)\in \mathbb{R}^l\times \mathbb{R}^l$ for $x\in
S^{n+1}(1)$. The condition $\langle Px, x\rangle=1$ is equivalent to
\begin{equation*}
|z|^2-|w|^2=1,
\end{equation*}
which implies $|z|^2=1,~|w|^2=0$. On the other hand, we observe that
\begin{eqnarray*}
V_+&:=&\{x\in M_-~|~\langle Px,
x\rangle=1\}\\
&=&\{x\in S^{2l-1}~|~\langle Px, x\rangle=1\}.
\end{eqnarray*}
Thus we get an isometry
\begin{equation*}
V_+\underset{isom.}{\cong}S^{l-1}(1).
\end{equation*}
Similarly, \begin{equation*} V_-\underset{isom.}{\cong}S^{l-1}(1).
\end{equation*}

Now the proof of Theorem \ref{thm2} is complete.

\end{document}